\documentclass[a4paper,11pt]{article}
\usepackage{epsfig}
\usepackage{amsmath}
\usepackage{amssymb}
\usepackage{amscd}

\newcommand{\newnotion}[1]{{\bf \emph{#1}}\index{#1}}

\newtheorem{theorem}{Theorem}[section]

\newtheorem{lemma}[theorem]{Lemma}

\newtheorem{definition}[theorem]{Definition}

\newtheorem{example}[theorem]{Example}
\newtheorem{commentary}[theorem]{Commentary}

\newtheorem{remark}[theorem]{Remark}

\begin{document}
\title{Veronese curves and webs Interpolation}
\author{T.B. Bouetou \& J.P. Dufour}

\maketitle

\begin{abstract}
 In this paper, we review basic results,
essentially due to J. Turiel, concerning the link between
classical multidimensional webs and Veronese webs.
\end{abstract}

\begin{quote}{\bf Keywords:}  Classical webs, (p+1)-webs,
Veronese curves, Veronese webs, bihamiltonian systems
\end{quote}
\begin{quote}{\bf AMS subject classification 2000:}
  53A20, 53C12, 53D30, 53D99.
\end{quote}

\section{Introduction.}
  In this survey, we will be interested in Veronese webs (particular case of  one
  parameter families of foliations), as
defined by \cite{gz,tu}, and  ordinary webs (finite families of
foliations in general position), as defined by Blaschke, Akivis
and Goldberg \cite{bla1,bla2,ak1,ak2,gol}. If we look  the
literature about webs, these two domains  were developed
apparently independently. Our main goal is to establish the link
between the two domains. \subsection{ Classical webs.}
\begin{definition}
  A k-web of codimension $c$ over a manifold $V$ is a
family of $k$ foliations $\mathfrak{F}_1 , \mathfrak{F}_2 ,\dots
,\mathfrak{F}_k  $ of $V,$ all of codimension $c,$ in ``general
position''. This last condition
 means that, if we denote by  $\mathfrak{F}_i(m)$ the tangent  plane
to $\mathfrak{F}_i$ at the point $m$ (the contact element to
$\mathfrak{F}_i$ at $m$), then the $c-$codimensional subspaces
  $\mathfrak{F}_1(m) , \mathfrak{F}_2(m) ,\dots ,\mathfrak{F}_k(m)  $ of
$\ T_{m}V$ are in general position  (as transverse as possible).
\end{definition}

  Classical examples are k-webs of $\mathbb{R}^2 ,$ i.e.  systems of $k$
 families of curves in a plane.
 We will call simply \newnotion{3-web}
any 3-web of  codimension $c$ on a $2c$-dimensional manifold.
General position means  here that we have $\mathfrak{F}_{i}(m)
\cap \mathfrak{F}_{j}(m)=\{ 0\} $ for $i \neq j$ and at any  $m$.
These webs are related to binary laws: if $(x,y)\longmapsto  x
\circ y $ is a smooth binary law on the manifold  $M,$ then we can
associate the three foliations $x=C^{st}, y=C^{st}$ et $x \circ
y=C^{st}$
 on $M \times M $; with further  hypotheses, this gives a 3-web (for example
if $x \circ y$ is a Lie  group law...). The 3-web were intensively
investigated in \cite{ak1,ak2,bla1,bla2,du,gri,na}.

 More generally a \newnotion{($p+1$)-web} will be a ($p+1$)-web of  codimension $c$ on
 a manifold of dimension $pc.$
   Here general position means that, for all $m \in V,$ we have
$$\mathfrak{F}_{1}(m) \cap \dots \cap \mathfrak{F}_{i-1}(m)\cap
\mathfrak{F}_{i+1}(m)\cap\dots \cap \mathfrak{F}_{p+1}(m)=\{
0\},$$ for any $i=1,\dots ,p+1.$
\begin{remark}
  A $p$-web of codimension $c$ on a manifold of dimension $pc$ is
locally  trivial, i.e. we can find  local coordinates
$$x^{1}_{1},\dots ,x^{1}_{c},\dots ,x^{p}_{1},\dots ,x^{p}_{c}$$
where $\mathfrak{F}_i$ is given by the equations $ \{
x^{i}_{1}=C^{te},\dots, x^{i}_{c}=C^{te} \}.$ So the first webs
which have an interesting local geometry are the above defined
($p+1$)-webs.
\end{remark}

\subsection{Veronese webs.}

\begin{definition}
   Let $V$ be a real  vector space of dimension
$(n+1)$. A  Veronese curve in the  projective space
$\mathbb{P}(V)$ is  a map
$$\gamma  : \mathbb{P}^1\mathbb{R} \longrightarrow \mathbb{P}V$$ which is the
     quotient of a map  of the type
$$
(x,y) \longmapsto x^n v_n+x^{n-1} yv_{n-1} +\dots +y^n v_0
$$
where $(v_0, v_1, \dots ,v_n)$ is a  base in $V$.
\end{definition}

\begin{definition} \label{veroneseweb}
 \cite{gz,pa2,tu1} A Veronese web of codimension
$c$ on a manifold $V$ of dimension $pc$ is a  one parameter family
of foliations $(\mathfrak{F}_t)_{t \in \mathbb{P}^{1}\mathbb{R}}$
of codimension $c$ on $V$ such that, for all $m$ in $V$, the
contact element $\mathfrak{F}_{t}(m)$ is given by:
$$\alpha^{1}_{t}=0,\dots, \alpha^{c}_{t}=0$$ with
$$\alpha^{i}_{t}=\gamma^{i}_{0}+t\gamma^{i}_{1}+t^2 \gamma^{i}_{2}+\dots +t^{p-1}\gamma^{i}_{p-1}$$
where $(\gamma^{i}_{j})_{^{i=1,\dots ,c}_{j=0,\dots ,p-1}}$ form a local
 coframe; that is
$$\gamma^{1}_{0}, \dots ,\gamma^{1}_{p-1}, \gamma^{2}_{0},\dots,
\gamma^{2}_{p-1},\dots ,\gamma^{c}_{0},\dots , \gamma^{c}_{p-1}$$
 are  differential forms, defined in a neighborhood of $m$ such
that $$\gamma^{1}_{0}(m), \dots ,\gamma^{1}_{p-1}(m),\dots
,\gamma^{c}_{0}(m),\dots , \gamma^{c}_{p-1}(m)$$ is a basis of
$T^{\star}_{m}V$.
\end{definition}
 Gelfand and Zakharevich \cite{gz} defined Veronese webs of  codim
$1$ $(c=1)$ and the notion was generalized by A. Panasyuk and J.
Turiel. In the sequel, we will sketch the way these notions
appeared.

 In the   bihamiltonian ``mechanic'', we study
pencil of Poisson structures  $\Pi_{0} +t\Pi_{\infty}$ over a
manifold
 $W$.  This means that  (\cite{gz})   $\Pi_{t}=\Pi_{0} +t\Pi_{\infty}$ are Poisson structures for all
 $t$ and  $ \Pi_{\infty}$ is also a Poisson structure.  It is equivalent to say
 that  $\Pi_{0}$ and $\Pi_{\infty}$ are Poisson structures with
 $[ \Pi_{0} ,\Pi_{\infty}]=0$ ($[\cdot ,\cdot ]$ is the so called Schouten bracket);
  in that case we say that  $\Pi_{0}$ and $ \Pi_{\infty}$ are ``compatible''.

For some time it was believed that any integrable Hamiltonian
system was a bihamiltonian system, i.e.  that  there exist a
second Poisson structure compatible with the Poisson structure
related to the initial symplectic structure, which should be
invariant by the Hamiltonian field. The correct idea is that, any
bihamiltonian system  is integrable but  R.Brouzet \cite{brou} has
shown that the former belief was wrong. Nevertheless the classical
integrable systems are all bihamiltonians.

 J. Turiel \cite{tu3}
has classified the pairs of compatible  Poisson structures
 $( \Pi_{0} ,\Pi_{\infty})$ with $ \Pi_{0} $  symplectic (here we are in a
even dimensional situation). On the other hand I. Gelfand and I.
Zakharevich were the first to investigate the odd dimensional
case. Precisely they consider a  pencil
$$\Pi_{t} =\Pi_{0}+t \Pi_{\infty}$$ on a $2p-1$-dimensional
manifold such that $\Pi_{t}$ is, for all $t,$ of  maximum rank
$(2p-2)$. The symplectic foliation  $\mathfrak{F}_{t}$ of
$\Pi_{t}$ is then of codimension $1$ and locally given by the
zeroes  of a form $\alpha_t$. It is not yet a Veronese foliation
in the sense of definition \ref{veroneseweb}, but we will explain
hereafter that it is the case up to a quotient.

In fact we have the following local models
$$\Pi_{0}(m)=e_{1} \wedge f_{1}+e_{2} \wedge f_{2}+\dots +e_{p-1} \wedge f_{p-1}$$
$$\Pi_{\infty}(m)=f_{1} \wedge e_{1}+f_{2} \wedge e_{3}+\dots +f_{p-1} \wedge e_{p}$$
where $ e_{1},\dots , e_{p}, f_{1},\dots , f_{p-1}$ is a well chosen  base
 of  $T_{m}V$;
denote by\\ $ e^{*}_{1},\dots , e^{*}_{p}, f^{*}_{1},\dots , f^{*}_{p-1}$
 the dual  base of  $T^{*}_{m}V$.
 The distribution  $\mathfrak{F}_{t}(m)$ is the symplectic foliation of
$$\Pi_{t}(m)=e_{1} \wedge f_{1}+e_{2} \wedge f_{2}+\dots +e_{p-1}
\wedge f_{p-1}+t(f_{1} \wedge e_{1}+f_{2} \wedge e_{3}+\dots +f_{p-1} \wedge e_{p}).$$
It is easy to see that the  distribution annihilates the form
$$\beta_{t}=e^{*}_{p}+te^{*}_{p-1}+\dots +t^{p-1}e^{*}_{1}.$$
and that   $\mathfrak{F}_{t}(m)$ contains $<f_{1},\dots
,f_{p-1}>.$


Take a submanifold  $V$ of dimension $p$ transverse to
$<f_{1},\dots ,f_{p-1}>$, the traces of  $\mathfrak{F}_{t}$ on $V$
form a Veronese web of codim $1$ defined by:
$$\alpha_{t}=e^{*}_{p}+te^{*}_{p-1}+\dots +t^{p-1}e^{*}_{1}.$$

The theory initiated by  Gelfand-Zakharevitch and ended by J.
Turiel says that the
 local  invariants of the pair $( \Pi_{0} ,\Pi_{\infty})$ are the local invariants of
 this Veronese foliation restricted to  $V$.
Latter the pairs   $( \Pi_{0} ,\Pi_{\infty})$  such that
 $ \Pi_{t}$ is of  constant corank $c> 1$ where investigated and, by the use of the same method, one
obtain Veronese webs in the sense of  definition
\ref{veroneseweb}.

\section{Link between $(p+1)$-webs and Veronese webs.}

Let    $(\mathfrak{F}_{t})_{t}$ be a Veronese web of codimension
$c$ over the ${pc}$-dimensional manifold $V.$ Assume that
$t_{1},\dots ,t_{p+1}$
 are  two by two distinct  then  $(\mathfrak{F}_{t_i})_{i=1,\dots ,p+1}$ gives a $(p+1)$-web:
In fact  $(\mathfrak{F}_{t_i})$ is locally given by
$$ \alpha^{1}_{t_i}=0,\dots , \alpha^{c}_{t_i}=0, $$
with
$$\alpha^{i}_{t}=\gamma^{i}_{0}+\gamma^{i}_{1}t+\dots +\gamma^{i}_{p-1}t^{p-1}.$$
Since
$$
\vmatrix
1 & t_{i_{1}} & \ldots & t_{i_{1}}^{p-1}\\
1 & t_{i_{2}} &  \ldots & t_{i_{2}}^{p-1}\\
\vdots & \vdots & \vdots & \vdots \\
1 & t_{i_{p}} & \ldots & t^{p-1}_{i_{p}}\\
\endvmatrix
$$
is a  Van Der Monde determinant, it is  clear that:
$$ \alpha^{1}_{t_{i_{1}}},\dots, \alpha^{c}_{t_{i_{1}}},\dots ,
 \alpha^{1}_{t_{i_p}},\dots, \alpha^{c}_{t_{i_p}} $$
form a  base of $T^{*}V$ for all $i_{1},\dots ,i_{p}$,  two by two
distinct with
 $\big\{i_{1},\dots ,i_{p}\big\} \subset \big\{1,\dots,p+1\big\} $,
 therefore we have the condition of general  position.

The most difficult problem is the passage from the $(p+1)$-webs to
Veronese webs. We have a problem of \newnotion{interpolation} of
     $(\mathfrak{F}_{i})_{i=1,\dots,p+1}$  to a curve
 $(\mathfrak{F}_{t})_{t \in \mathbb{P}^{1}\mathbb{R}}$  having good properties.
We decompose this into two  problems. \begin{itemize}

 \item {\bf Algebraic interpolation:}
Given $p+1$ subspaces of codimension $c$ in a $pc$-dimensional
vector space $V$ in general position, find a natural curve
  of subspaces of codimension $c$ in  $V$  passing through the given
 $p+1$ subspaces. It is a problem in $G_{c}( V)$ the
 Grassmannian of subspaces of codimension $c$ in $V$. Let us
assume that this problem has a unique solution.  Given $p+1$
distributions of
 codimension $c$ on a manifold of dimension $pc$  there would exist a
natural method to interpolate these $p+1$ distributions
 $\mathfrak{F}_{1},\mathfrak{F}_{2},\dots ,\mathfrak{F}_{p+1}$ in a curve
$\mathfrak{F}_{t}$ of distributions.

\item{\bf Integrability:} Under which condition these
distributions are integrable? For example, is the integrability of
$\mathfrak{F}_{1},\mathfrak{F}_{2},\dots,\mathfrak{F}_{p+1}$
sufficient to guaranty that of $\mathfrak{F}_{t}$ for all
$t$?\end{itemize}

In the next sections we will  investigate these questions.

\section{Interpolation of a finite family of  subspaces.}

Let $V$ be a vector space and $G_{c}( V)$ the Grassmannian of his
codimension $c$ subspaces. We put $N=$dim$V-c$ and denote by
$S_c(V)$ the open subset of $V^N$ formed by $N$-uples of linearly
independent vectors of $V.$ Let $\beta : \mathbb{P}^{1}\mathbb{R}
\longrightarrow G_{c}( V)$; it is said to be a degree $q$ curve if
it pulls back as:
 $$\hat{\beta} : \mathbb{R}^{2}\setminus 0  \longrightarrow S_c(V)$$
$$\hat{\beta}(x,y)=\big( \beta_{1}(x,y),\dots , \beta_{N}(x,y)  \big),$$
where $\beta_i$ has the form
$$ \beta_{i}=\sum_{j}\beta^{j}_{i}(x,y)e_j,$$
$(e_j)_j$ is a basis of $V$ and $\beta^{j}_{i}$ are homogeneous
polynomials of degree $q$; we have the following commutative
diagram

$$
\begin{CD}
\mathbb{R}^2\setminus {0} @>{\hat{\beta}}>>S_c(V)\\
@V{P}VV  @VV{P}V \\
\mathbb{P}^1 \mathbb{R}   @>{\beta}>> G_{c}(V)
\end{CD}
$$
where $P$ are canonical projections $\big(P(v_{1},
\dots,v_{n})=<v_{1}, \dots,v_{n}>\big)$.

 Let $F_1,\dots ,F_{p+1}$ be  given points of $G_{c}(V);$ we will say that
$\beta : \mathbb{P}^{1}\mathbb{R} \longrightarrow G_{c}( V)$ is a
\newnotion{minimal interpolation} of  $(F_1,\dots ,F_{p+1})$ if $\beta $ is
a curve of minimal  degree $q$  passing through   $F_1,\dots
,F_{p+1}$.

 It is a difficult problem to find such minimal interpolation and
see if they are unique: in general it is wrong. Furthermore these
curves are not independent of the choice of the
\newnotion{parametrization}: the sequences of $t_i$ such that
$\beta(t_i)=F_i$. In the sequel we will show that there are unique
minimal interpolations in two important cases:
\begin{itemize} \item the case where $\dim V=pc$ ($\dim
F_i=(p-1)c$)
 \item the case where $\dim V=pN$ ($\dim F_i=N$).
 \end{itemize}
In the first case  the minimal interpolations are
\newnotion{pencils}, i.e. degree 1
curves, of $c$-codimensional subspaces; in the second case we
recover Veronese curves and their generalization. Moreover  these
cases are dual to each other.

\subsection{Pencils interpolations.}

In this paragraph we deal with a family $F_1,\dots ,F_{p+1}$ of
$c$-codimensional subspaces of the $pc$-dimensional vector space
$V.$ We suppose that this family is in general position: this
means that, for every $i,$ we have $F_1\cap\dots\cap F_{i-1}\cap
F_{i+1}\cap\dots\cap F_{p+1}=\{ 0\} .$

Fix a system of linear coordinates $(x_1^1,\dots
,x_1^c;x_2^1,\dots ,x_2^c;\dots ;x_p^1,\dots ,x_p^c)$ such that
the equations of $F_i$ are $x_i^1=0,\dots ,x_i^c=0$ for $i=1,\dots
,p,$ and $F_{p+1}$ has equations $\sum_ix_i^1=0,\dots
,\sum_ix_i^c=0.$ We denote by
$$(e_1^1,\dots ,e_1^c;e_2^1,\dots
,e_2^c;\dots ;e_p^1,\dots ,e_p^c)$$
the corresponding basis.

We fix also a system $t_1,\dots ,t_p$  of two by two different
real numbers.

A pencil of $c$-dimensional subspaces is a degree 1 curve $\beta$
of $c$-dimensional subspace of $V$: this means that $\beta$ pulls
back as a curve $\hat\beta :\mathbb{R}^{2}\setminus 0 \rightarrow
S_c(V)$ with $\hat\beta (x,y)= (xa^1+yb^1,\dots
,xa^{(p-1)c}+yb^{(p-1)c}),$ where $a^j$ and $b^j$ are vectors of
$V.$ With the identification $t\equiv [t:1],$ we can write $\beta
(t) =(G-t{\rm Id})\beta (\infty ),$ where $\beta (\infty )$ is the
space generated by the $a^j$ and $G:V\rightarrow V$ is any linear
map such that $G(a^j)=-b^j$ for every $j.$

 We want to interpolate the $F_i$ by such a pencil. More precisely,
 we want a pencil $\beta$ with $\beta (t_i )=F_{i},$ for $i=1,\dots ,p$
 and $\beta (\infty )=F_{p+1}.$ A simple solution is obtained by
 choosing $G$ such that $G(e^j_i)=t_ie^j_i$ for every $i=1,\dots
 ,p$ and $j=1,\dots ,c:$ we have
 $$\beta (\infty )=<e^j_p-e^j_k;k=1,\dots ,p-1;j=1,\dots ,c>,$$
 then
 \begin{equation}\label{betat}
\beta (t )=<(t_p-t)e^j_p-(t_k-t)e^j_k;k=1,\dots ,p-1;j=1,\dots
,c>,
\end{equation}
and it is easy to see that $\beta (t_i)$ has equations
$x_i^1=0,\dots ,x_i^c=0.$

In the sequel we will investigate the uniqueness of this pencil.

First we will suppose there is another linear map $G'$ with
$(G'-t_i{\rm Id})(F_{p+1})=F_{i}$ for every $i=1,\dots ,p.$ Put
$$G'(e^j_p-e^j_k):=u^j_k=\sum_{r=1\dots c,s=1\dots
p}a^{jr}_{ks}e^r_s.$$ Then equations
$x^j_r(u_i^s-t_r(e_p^s-e_i^s))=0$ for every $s,j=1,\dots ,c,$
$i=1,\dots ,p-1$ and $r=1,\dots ,p,$ lead to
$$u^j_k=t_pe^j_p-t_ie^j_k$$
for every $j=1,\dots ,c$ and $k=1,\dots ,p-1.$ So the pencil
attached to $G'$ is exactly $\beta $ (the one attached to $G$).

We can remark that the difference $\Delta =G'-G$ is a linear
mapping of $V$ such that $\Delta (e^j_p-e^j_k)=0$ for every $j$
and $k.$ So $\Delta$ is characterized by the fact that there are
arbitrary vectors $v^1,\cdots ,v^c$ of $V$ with $\Delta
(e^j_k)=v^j,$ for every $j$ and $k.$ In particular we can always
manage such that $G'$ is invertible: if the $t_i$ are all non
zero, then $G$ is invertible; if, for example, $t_1$ vanishes, we
can choose $v^j=e^j_1$ for every $j.$

Next we remark that coordinates $(x_i^j)^{j=1,\dots ,c}_{i=1,\dots
,p}$ are unique up to a linear change of the form
$(x_i^j)'=\sum_{s=1,\dots ,c}a_s^jx_i^s;$  this means that the
matrix of this linear change is a $pr\times pr$ matrix which have
only null terms except $p$ diagonal $r\times r$
  blocs all equal to $A=(a^r_s)_{r,s=1,\dots ,p}.$ This
 induces that $\beta$ doesn't depend on the particular choice of
 the adapted coordinates $x_i^j.$

 Finally we want to see how this interpolation depends on the
 parametrization, i.e. on the sequence $t_1,\dots ,t_p.$ First of
 all, remark that, if $\beta$ is a pencil as above, then we can
 perform a projective transform on the parameter space
 $\mathbb{P}^{1}\mathbb{R}$ and we keep a pencil. This allows us
 to impose the values at three different points: this justifies a
 posteriori the particular choice of $\beta (\infty )$ in the
 preceding calculations. We could also have fixed two other
 values, for example $t_1=0$ and $t_2=1$ (imposing $\beta (0)=F_1$ and $\beta
 (1)=F_2$). The following lemma says that two pencils which
 interpolate $F_1,\dots ,F_{p+1}$ are the same if and only if the
 sequences of parameters $\tau_1,\dots ,\tau_{p+1},$ where these
 pencils pass respectively at $F_1,\dots ,F_{p+1},$ are the same
 up to a projective
 transformation of $\mathbb{P}^{1}\mathbb{R}.$

 \begin{lemma} \label{unicity} Let $\beta$ and $\beta'$ be two pencils,
 interpolating $F_1,\dots ,F_{p+1},$ such that
 $$\beta (\infty )=\beta'(\infty )=F_{p+1},\ \
 \beta (0 )=\beta'(0)=F_{1},\ \ \beta (1 )=\beta'(1 )=F_{2}.$$
 Let $t_i$ and $t'_i$ for $i=3,\dots
 ,p$ the values of the parameters such that $F_i=\beta (t_i)=\beta'(t'_i ).$
  Then  $\beta$ and $\beta'$ have the same image ($\{\beta (t );t\in
\mathbb{P}^{1}\mathbb{R}\}=\{\beta' (t );t\in
\mathbb{P}^{1}\mathbb{R}\}$) if and only if $t_i=t'_i$ for every
$i=3,\dots ,p.$\end{lemma}

\begin{Proof}
The preceding calculations give the ``if'' part. To prove the
converse we suppose that, for each $t\in
\mathbb{P}^{1}\mathbb{R},$ there is $t'\in
\mathbb{P}^{1}\mathbb{R},$ with
$$\beta (t)= \beta'(t').$$ Formula \ref{betat} gives
$$<(t_p-t)e^j_p-(t_k-t)e^j_k;k=1,\dots ,p-1;j=1,\dots
,c>=$$$$<(t'_p-t')e^j_p-(t'_k-t')e^j_k;k=1,\dots ,p-1;j=1,\dots
,c>,$$ for every $t.$ From this we deduce equations
$$(t_p-t)(t'_k-t')=(t'_p-t')(t_k-t),$$
 so  relations
$$t'=t\frac{t'_p-t'_k}{t_p-t_k}+\frac{t_pt'_k-t'_pt_k}{t_p-t_k},$$for $k=1,\dots ,p-1.$
Then hypothesis $t_1=t'_1$ and $t_2=t'_2$ imply $t_k=t'_k$ for
every $k.$
\end{Proof}

 \subsection{Veronese interpolations.}

To each subspace $F$ of the vector space $V$ we associate its
annihilator $F^\circ$ which is the subset of $V^*$ formed by the
linear forms on $V$ which vanish on $F.$ Now if $\beta $ is an
interpolation of the family of subspaces $F_1,\dots ,F_{p+1}$ of
$V,$ then $\beta^\circ,$ defined by $\beta^\circ(t)=(\beta
(t))^\circ, $ is an interpolation of the family of subspaces
$F_1^\circ,\dots ,F_{p+1}^\circ.$

Now suppose that $V$ has dimension $pN$ and the $F_i$ have
dimension $N.$ Then $F_i^\circ$ have codimension $c:=N$ and the
preceding subsection gives pencil interpolations, in the general
position cases, for $F_1^\circ,\dots ,F_{p+1}^\circ.$ Denote by
$\gamma$ such a pencil; we have
 (formula \ref{betat})
 $$\gamma (t )=<(t_p-t)\alpha^j_p-(t_k-t)\alpha^j_k;k=1,\dots ,p-1;j=1,\dots
,c>,$$ for a good  basis $(\alpha^j_k)_{i,k}$ of $V^*$ and a
parametrization such that
$$\gamma (\infty )=F_{p+1}^\circ,\ \ \gamma (t_i )=F_i^\circ$$
for $i=1,\dots ,p.$ If $(e^j_k)_{i,k}$ is the dual basis to
$(\alpha^j_k)_{i,k},$ we have
$$\gamma^\circ (t )=<\sum_{i=1}^p(t_1-t)\cdots(t_{i-1}-t)(t_{i+1}-t)\cdots(t_1-t)e^j_i;j=1,\dots
,c>.$$

So we get a degree $p-1$ interpolation of $F_1,\dots ,F_{p+1}.$ In
the case where $N(=c)=1,$ we can prove that $\gamma^\circ$ gives a
Veronese curve in $\mathbb{P}(V).$ For this reason we call these
$\gamma^\circ$ \newnotion{Veronese interpolations}, even in the
case $N>1.$ Uniqueness properties of pencil interpolations
translate into corresponding uniqueness properties for Veronese
interpolations.

\begin{example} For $p=2$ a Veronese interpolation is also a pencil.
\end{example}
\begin{example} For $p=3$ and $N=1$ a Veronese interpolation is a degree 2 curve
in a projective plane: it is a conic. We recover that there are
conics passing by four given points, and the lemma \ref{unicity}
is the generalization of the classical result which says that such
a conic is characterized by the cross-ratio of these four points
on the conic.
\end{example}

\section{Integrability of distributions.}

\subsection{Distributions.}

 The results of the preceding section are purely  algebraic but they
  pass to smooth distributions on manifolds.
 For example, when we have $(p+1)$ smooth distributions
$\mathfrak{F}_1 ,\dots ,\mathfrak{F}_{p+1}$ of codimension $ c,$
in general position, on a \newnotion{manifold} $W$ of dimension
 $ pc$, we can work point by point in each tangent space  $T_{m}W$ to construct the distribution
$\mathfrak{F}_{t}$ which interpolate them. The uniqueness of this
procedure ensures their smoothness. To be coherent with the
vocabulary of our second section, we call these 1-parameter
families of distribution \newnotion{Veronese distributions}.

 In the
neighborhood $\mathcal{U}$ of each point $m$ we have a family of
operators $G(m),$ depending smoothly on $m,$ such that
$$\mathfrak{F}_{t}(m)=(G(m)-tI)\mathfrak{F}_{\infty}(m).$$

\subsection{Integrability theorem.}
  In this  subsection, we will give a short proof of the following theorem of A.
  Panasiuk (\cite{pa1}).

\begin{theorem}
 Let $\big(\mathfrak{F}_{t} \big)_{t}$ be a Veronese
distribution on a $pc$-dimensional manifold  $W$. The distribution
$\mathfrak{F}_{t}$
 is integrable for any $t$ if and only if there exist $p+2$ values of $t$ for
which   $\mathfrak{F}_{t}$ is integrable.
\end{theorem}

 This theorem is not evident in the covariant  version, i.e. when we
 define distributions as zeroes of set of forms $\alpha^{1}(t),\dots,
 \alpha^{c}(t):$
using the Frobenius theorem, the integrability of
$\mathfrak{F}_{t},$ for any $t,$ is locally equivalent to
$$d \; \alpha^{i}(t)\wedge \alpha^{1}(t)\wedge ...\wedge \alpha^{c}(t) \equiv 0 ,$$
for all $t.$ This gives a polynomial equation of  degree $(c+1)p$
in $t.$ It will vanish identically if it vanishes at $(c+1)p+1$
values of $t$ which is, in general, bigger than $p+2.$

It is not also evident in  contravariant version, i.e. when we
define distribution by means of vector fields: $\mathfrak{F}_{t}$
is integrable, for any
  $t,$  if and only if, for any $t,$
$$\big[ X_{i}(t), X_{j}(t)\big]\wedge  X_{1}(t)\wedge ... \wedge X_{(p-1)c}(t)=0,$$
by denoting $\mathfrak{F}_{t}=<X_{1},...,X_{(p-1)c}>$ where
$X_{1},...,X_{(p-1)c} $ form a local basis. This gives a
polynomial equation of degree $2+(p-1)c$. It will vanish
identically if it vanishes at   $3+(p-1)c$ values of $t,$ still
bigger than $p+2.$

\begin{Proof}
It is sufficient to work locally in a neighborhood of any point of
$W$: we choose invertible operators $G(m),$ depending smoothly on
$m,$ with $\mathfrak{F}_{t}(m)=(G(m)-tI)\mathfrak{F}_{\infty}(m).$
 We choose also a family of vector fields $ \big(
v_{1},\dots,v_{(p-1)c}\big)$ which generates locally
$\mathfrak{F}_{\infty}(m).$

The integrability of $\mathfrak{F}_{t}$ is given by  equation:

$$\big[ (G-tI)v_{i},  (G-tI)v_{j}\big]=\sum_{k}\theta_{ij}^{k} (G-tI)v_{k}$$
for any $ i,j$.

Let assume that $\big( \mathfrak{F}_{t}\big)$ is integrable for
$p+2$ values of $t$; we can  assume that it is true for $t=0,
\infty$ et $t_{1},...,t_{p}$ two by two  distinct. This implies
relations
$$ \big[ v_{i}, v_{j}\big]=\sum_{k}\alpha_{ij}^{k}v_{k}$$
(for $t=\infty$),
$$ \big[G v_{i},G v_{j}\big]=\sum_{k}\beta_{ij}^{k}G\big(v_{k}\big)$$
(for t=0). We have
\begin{equation}\label{integ}
 \big[ (G-tI)v_{i},
(G-tI)v_{j}\big]=\big[Gv_{i},Gv_{j}\big]-t\Delta
\big(v_{i},v_{j}\big)+t^{2}\big[v_{i},v_{j}\big],
\end{equation}
with $ \Delta
\big(v_{i},v_{j}\big)=\big[Gv_{i},v_{j}\big]+\big[v_{i},Gv_{j}\big]$.
We introduce the Nijenhuis torsion $N_{G}$ \cite{nij} of $G$:
$$ N_{G}\big( v_{i},v_{j}\big)=\big[Gv_{i},Gv_{j}\big]-
G\Delta \big(v_{i},v_{j}\big)+G^{2}\big[v_{i},v_{j}\big].$$ Then
we get
 $$ \Delta \big(v_{i},v_{j}\big)=G^{-1}\big[Gv_{i},Gv_{j}\big]+G\big[v_{i},v_{j}\big]-G^{-1}N_{G}\big( v_{i},v_{j}\big).$$
Thus the first member of formula  \ref{integ} becomes
$$ (I-tG^{-1})\big[Gv_{i},Gv_{j}\big]-  t\big((G-tI)\big[v_{i},v_{j}\big]\big)+tG^{-1}N_{G}\big( v_{i},v_{j}\big)$$
$$ =G^{-1}(G-tI)\big[Gv_{i},Gv_{j}\big]-t\big((G-tI)\big[v_{i},v_{j}\big]\big)+tG^{-1}N_{G}\big( v_{i},v_{j}\big)$$
$$ =G^{-1}(G-tI)\sum_{k} \beta_{ij}^{k}G\big(v_{k}\big)- t\bigg((G-tI)\sum_{k} \alpha_{ij}^{k}v_{k}\bigg)+tG^{-1}N_{G}\big( v_{i},v_{j}\big)$$
$$ =(G-tI)\sum_{k} \gamma_{ij}^{k}(t)v_{k} +tG^{-1}N_{JG}\big( v_{i},v_{j}\big)$$
with $\gamma_{ij}^{k}(t)= \beta_{ij}^{k}-t\alpha_{ij}^{k}.$ The
integrability for $t=t_{1}, t_{2},\dots ,t_{p}$ gives us equations
$$ t_{r}G^{-1}N_{G}\big( v_{i},v_{j}\big)=\big( G-t_{r}I\big)\big(\sum_k \mu_{ij}^{k}\big(t_{r} \big) v_k\big)
$$
for $r=1,....,p$.\\
Therefore  $G^{-1}N_{G}\big( v_{i},v_{j}\big)$ is in $
\bigcap_{r=1}^{p}\mathfrak{F}_{t_{r}} .$ As we have
$$\bigcap_{r=1}^{p}\mathfrak{F}_{t_{r}}=\{0\}, $$
($t_{r}$ two by two  distinct), we can conclude

$$ G^{-1}N_{G}\big( v_{i},v_{j}\big)=0, $$
then  $N_{G}\big( v_{i},v_{j}\big)=0,$  for any $ i,j$. So we
have, for any  $t$,

$$\big[ (G-tI)v_{i},  (G-tI)v_{j}\big]=\sum_{k}\theta_{ij}^{k} (G-tI)v_{k}.$$
So we obtain the integrability of each $\mathfrak{F}_{t}$.
\end{Proof}

\begin{remark} In his study of Veronese webs (see \cite{tu2,tu})
J. Turiel invented the above technics. He fixes $p+1$ foliations
of the family, say $\mathfrak{F}_{\infty}$ and
$\mathfrak{F}_{t_i},$ for $i=1,\dots ,p.$ The $p$ distributions
$\mathfrak{H}_{i},$ defined by
$$\mathfrak{H}_{i}(m)=\bigcap_{j=1,\dots,i-1,i+1,\dots ,p}\mathfrak{F}_{j}(m),$$
for $i=1\dots ,p,$   are integrable and decompose, at each point
$m,$ the tangent space in a direct sum. So the operator $G$
defined by $G=t_iI$ in restriction to every $\mathfrak{H}_{i},$
has a null Nijenhuis torsion. This simplifies the above
calculations; nevertheless the integrability of
$\mathfrak{F}_{\infty}$ and $\mathfrak{F}_{t_i},$ for $i=1,\dots
,p$ doesn't  ensure that of the whole family because either $G$ is
not invertible or we don't know if $\mathfrak{F}_{0}$ is
integrable. As we saw in the paragraph preceding lemma
\ref{unicity}, we can replace our $G$ with $G'=G+\Delta$ for a
well chosen $\Delta ;$ doing this we can get $G'$ with non zero
Nijenhuis torsion.
\end{remark}


\addcontentsline{toc}{section}{References}
\begin{thebibliography}{99}

\bibitem{ak1}Akivis M.A., Goldberg V.V., Differential geometry of Webs,in Handbook of differential geometry edited by F.J.E. Dillen, L.C.A.Verstraelen, vol. 1, chapter 1, Elsevier Science B.V., (2000), pp. 1-152
\bibitem{ak2}  Akivis M.A., Shelekhov A.M. Geometry and algebra of multidimensional
3-Webs 1992, Kluwer academic publishers, Dordrecht-Boston-London, XVII+385p.

\bibitem{bla1} Blaschke W. Einfuhrung in die geometrie der waben (Introduction to webs goemetry) Birkhauser Verlag, Basel und Stuttgart 1955, 108p.

\bibitem{bla2} Blaschke W, Bol G. Geometrie der Gewebe, Grundlehren der Mathematischen Wissenschaften, Springer, 1938.

\bibitem{bol} Bol G \"Uber ein bemerkenswertes 5-Gewebe in der Ebene, Abh. Math. sem. hamburg {\bf 11} (1936), pp. 387-393.

\bibitem{brou} Brouzet R. About the existence of recursion operators for completely integrable hamiltonian system near Liouville torus, J. math. phys. {\bf 34} (1993) $n^{0} 4$ p. 1309-1313.

\bibitem{du} Dufour J.P. Introduction aux tissus, S\'eminaire Caston Darboux, Montpellier (1989/1990)

\bibitem{dou} Dubrovin B, ; Zhang Y. Normal forms of hierarchies of
integrable PDEs, Frobenius manifolds and Gromov-Witten invariants,
arXiv:math.DG/0108160 v1 23 august 2001.

\bibitem{fn} Fr\"olicher A.; Nijenhuis A. Theory of vector-valued differential forms I, Indag. Math., {\bf 18} (1956), 338-359.

\bibitem{gz} Gelfand I.M., Zakharevich I Webs, Veronese curves and bihamiltonian systems, J. Funct. Anal. {\bf 99} (1991), pp. 150-178.

\bibitem{gri} Grifone J., Salem E. Web theory and related topics, collected papers from the contributions to the Journ\'ee sur les Tissus held in Toulouse in December 1996, editet by J. Grifone and E. Salem World Scientific edition 2001. 234p. ISBN 981-02-4604-8.

\bibitem{gol} Goldberg V.V.Theory of multicodimensional (n+1)-Webs, Kluver,
Dordrecht, (1988)

\bibitem{na} Nagy P.T. Invariant tensorfields and the canonical connection of 3-web Aequationes Math. {\bf 35} (1988), pp.31-44

\bibitem{nij} Nijenhuis A. Jacobi-type identities for bilinear concomitants of certain tensor fields, Indagationes Maths.,{\bf A58} (1955), 390-403.
\bibitem{pa1} Panasyuk A. On integrability of generalized Veronese curves of distributions, Rep. Math. Phys. {\bf 50} (2002) $n^{0}$ 3, pp. 291-297.

\bibitem{pa2} Panasyuk A. Veronese webs for bihamiltonian structures of higher corank, Banach center pub. {\bf 51} pp. 251-261, Polish. Acad. Sci. Warsaw 2000.
\bibitem{mhr} Rigal M.H. G\'eometrie globale des syst\`emes bihamiltoniens en dimension inpaire Th\`ese Universit\'e de Montpellier

\bibitem{tu1} Turiel F.J., Tissus de Veronese analytiques de codimension sup\'erieure et structures bihamiltoniennes, C.R. Acad. Sci. Paris, t. 331, S\'erie I, p. 61-64, 2000.

\bibitem{tu2} Turiel F.J., $C^{\infty}$-\'equivalence entre tissus de tissus de Veronese  et structures bihamiltoniennes, C.R. Acad. Sci. Paris, t. 328, S\'erie I, p. 891-894, 1999.
\bibitem{tu} Turiel F.J., $C^{\infty}$-classification des germes de Veronese  et structures bihamiltoniennes, C.R. Acad. Sci. Paris, t. 329, S\'erie I, p. 425-428, 1999.
\bibitem{tu3} Turiel F.J Classification locale simultan\'ee de deux formes symplectiques compatibles, Manuscripta Math. {\bf 82} (1994), $n^0$ 3-4, p. 349-362.

\end{thebibliography}
\end{document}